# How to deal with numbers of decision making units and variables in data envelopment analysis

**Dariush Khezrimotlagh**[*]

Department of Applied Statistics, Faculty of Economics and Administration, University of Malaya, Kuala Lumpur, Malaysia, dariush@um.edu.my, January 01, 2014.

**Abstract**

Sufficient numbers of Decision Making Units (DMUs) in comparison with the number of input and output variables has been a concern of using Data Envelopment Analysis (DEA) in the last three decades. There are several studies in the literature of DEA which have tried to handle this issue by providing additional procedures to increase the number of DMUs, decreasing the number of variables or finding a relationship between the number of DMUs and variables. However, there are no concerns about the number of DMUs in comparison with the number of variables while Kourosh and Arash Method (KAM) is applied. A geometric reason is provided to depict the validity of the method without any extra conditions or additional methodologies. The technique is quite simple with no computational complexities of current methodologies even if the number of DMUs is less than the number of variables. A real-life numerical example of 32 DMUs with 45 variables demonstrates the advantages of the proposed technique.

**Keywords**: DEA, KAM, Efficiency, Technical efficiency, Ranking.

## 1. Introduction

Charnes *et al*. (1978) proposed Data Envelopment Analysis (DEA) to assess the relative efficiency of a set of homogenous Decision Making Units (DMUs) with multiple inputs and multiple outputs. DEA has been dramatically developed in the last three decades, and is known as a popular non-parametric technique which only requires a simple set of inputs and outputs values. A description on DEA's literature can be seen in Ray (2004).

Although, DEA provides a good number of new insights and additional information not available from conventional econometric methods (Seiford and Thrall, 1990), there is a concern about using DEA due to insufficient numbers of DMUs in comparison with the number of variables. Indeed, as the number of DMUs decreases or the number of variables increases, the technically efficient DMUs are increased which decreases the discriminative power of DEA.

There are several studies in the literature of DEA which have provided methodologies to reduce the number of variables, or find a relation between the number of DMUs and variables such as the studies of Golany and Roll (1989), Banker *et al.* (1989), Pedraja-Chaparro *et al*. (1999), Dyson *et al*. (2001), Jenkins and Anderson (2003), Cooper *et al*.

---

[*] e-mail address: khezrimotlagh@gmail.com

(2007), Morita and Avkiran (2009) and Osman *et al*. (2011). However, this study, unlike the previous studies, illustrates there is no concern about the numbers of DMUs in comparison with the number of variables while the new proposed technique is applied. The new technique can be measured by Kourosh and Arash Model (KAM) which was recently proposed by Khezrimotlagh *et al*. (2013a) to improve the DEA capabilities to distinguish between DMUs appropriately.

The rest of this paper is organized in five sections. Section 2 is a short background on KAM. Section 3 represents a methodology to distinguish between technically efficient DMUs without concern about the number of DMUs and variables. A real numerical example proposed by Osman *et al*. (2011) is represented in Section 4 to depict the advantages of KAM to deal with this issue and the paper is concluded in the last section. Simulations are also performed using Microsoft Excel Solver as it required simple linear programming.

## 2. Background on KAM

Suppose that there are $n$ DMUs ($\text{DMU}_i, i = 1,2,\ldots,n$) with $m$ non-negative inputs ($x_{ij}, j = 1,2,\ldots,m$) and $p$ non-negative outputs ($y_{ik}, k = 1,2,\ldots,p$), such that, at least one of the inputs and one of the outputs of each DMU are not zero, and for every $i$ there is a $j$ such that $x_{ij} \neq 0$ and also for every $i$ there is a $k$ such that $y_{ik} \neq 0$. The linear KAM is as follows for an appropriate $\epsilon = (\varepsilon^-, \varepsilon^+) \in \mathbb{R}_+^{m+p}$, where $\varepsilon^-$ is $(\varepsilon_1^-, \varepsilon_2^-, \ldots, \varepsilon_m^-)$ and $\varepsilon^+$ is $(\varepsilon_1^+, \varepsilon_2^+, \ldots, \varepsilon_p^+)$ (Khezrimotlagh *et al*. 2013a):

$$\max \ \sum_{j=1}^m w_j^- s_{lj}^- + \sum_{k=1}^p w_k^+ s_{lk}^+,$$
Such that
$$\sum_{i=1}^n \lambda_i x_{ij} + s_{lj}^- = x_{lj} + \varepsilon_j^-, \ \forall j, \quad \sum_{i=1}^n \lambda_i y_{ik} - s_{lk}^+ = y_{lk} - \varepsilon_k^+, \forall k,$$
$$x_{lj} - s_{lj}^- \geq 0, \ \forall j, \qquad y_{lk} + s_{lk}^+ - 2\varepsilon_k^+ \geq 0, \forall k,$$
$$\sum_{i=1}^n \lambda_i = 1; \qquad \lambda_i \geq 0, \ \forall i,$$
$$s_{lj}^- \geq 0, \ \forall j, \qquad s_{lk}^+ \geq 0, \ \forall k.$$

The best technical target and score with $\epsilon$-Degree of Freedom (DF) are respectively depicted as:

$$\begin{cases} x_{lj}^* = x_{lj} - s_{lj}^{-*} + \varepsilon_j^-, \forall j, \\ y_{lk}^* = y_{lk} + s_{lk}^{+*} - \varepsilon_k^+, \forall p, \end{cases} \qquad KA_\epsilon^{*l} = \frac{\sum_{k=1}^p w_k^+ y_{lk} / \sum_{j=1}^m w_j^- x_{lj}}{\sum_{k=1}^p w_k^+ y_{lk}^* / \sum_{j=1}^m w_j^- x_{lj}^*}.$$

After the optimization, the following definition is able to identify the efficient DMUs regarding the goals of evaluation (Khezrimotlagh *et al*. 2013a).

**Definition:** A technical efficient DMU is KAM efficient with $\epsilon$-DF in inputs and outputs if $KA_0^* - KA_\epsilon^* \leq \delta$. Otherwise, it is inefficient with $\epsilon$-DF in inputs and outputs. The proposed amount for $\delta$ is '$10^{-1}\epsilon$' or '$\epsilon/(m+p)$'.

For more information about KAM and how to apply it, see Khezrimotlagh *et al*. (2012a-f, 2013a-d) and Khezrimotlagh (2014a-b).

## 3. Deal with umbers of DMUs and variables

It has been consistently suggested in the literature of DEA that there should be sufficient numbers of observations in comparison with the numbers of factors. For example, Pedraja-Chaparro *et al*. (1999) said DEA loses its discrimination power in terms of number of technical efficient and inefficient units when the value of $n/(m + p)$ is too small. Golany and Roll (1989) suggested that $n$ should be greater than **2 × ($m + p$)**, whereas Banker *et al*. (1989), Friedman and Sinuany-Stern (1998) and Cooper *et al*. (2007) proposed that it should be greater than **3 × ($m + p$)**. Dyson (2001) recommended that $n$ should be greater than **2 × $m$ × $p$**. However, this section illustrates how KAM deals with this issue even if the numbers of DMUs is less than the numbers of variables, although, a good number of DMUs is more appropriate.

Suppose that there are two non-dominated DMUs A and B, while each one has a single constant input and two output values. Assume that the values of input and outputs are commensurate with the unity scale. Two different situations for A and B can be considered in output spaces, as depicted in Figures 1 and 2.

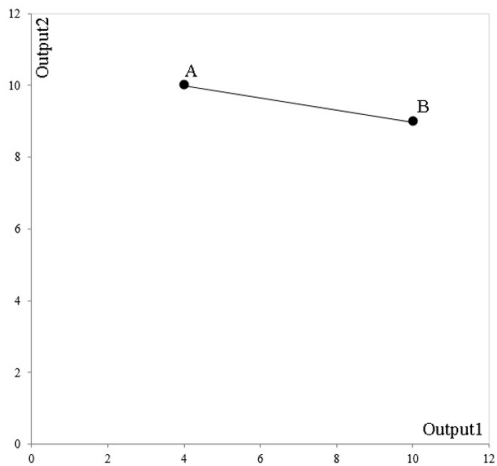
Figure 1: Two technically efficient DMUs (a).

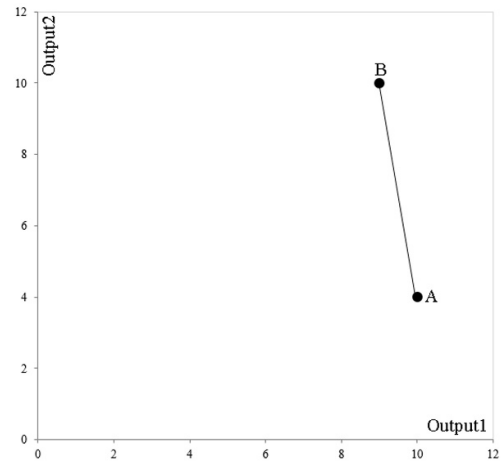
Figure 2: Two technically efficient DMUs (b).

Since the output values are considered commensurate with the unity scale, it is clear that the situation of B is better than the situation of A. In other words, B is more efficient than A, because the efficiency scores of A (in both figures) and B can be measured by (Output1+Output2)/Input, which are (4+10)/1 (=14) and (10+9)/1 (=19), respectively.

Now let's look at the neighbors of A in the Production Possibility Set (PPS), as depicted in Figures 3 and 4. In other words, assume that $A'$ is under evaluation where its components can be considered as (1; 3.5, 9.5) in Figure 3 and (1; 9.5, 3.5) in Figure 4.

Note that, $A'$ is a neighbor of A, meaning that the values of $A'$'s components have very small difference values with the values of A components. So, while very small errors are introduced in the components of A, such $A'$ can be selected. In this example, it is supposed that $\varepsilon_1^+ =$ **0.5** for the first output and $\varepsilon_2^+ =$ **0.5** for the second output. These values of epsilons are considered to have transparent figures and illustrate the method clearly. One may consider smaller values of epsilons and magnify the figures to make it clear. Indeed, this approach is independent of the values of epsilon, although, every value of epsilon has different meaning in conclusion.

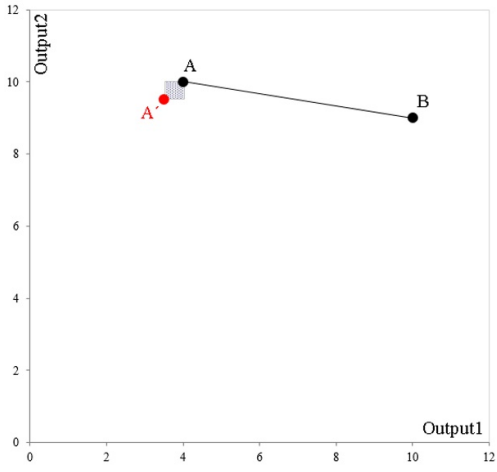
Figure 3: A neighbor of A (a).

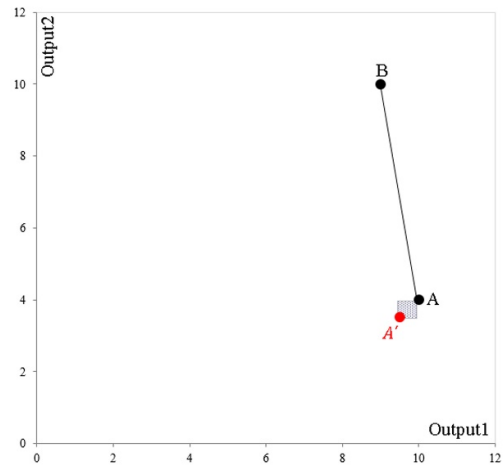
Figure 4: A neighbor of A (b).

As the arrows in Figures 5 and 6 depict, $A'$ does not suggest to A while optimum slacks are calculated, and it is suggested toward B.

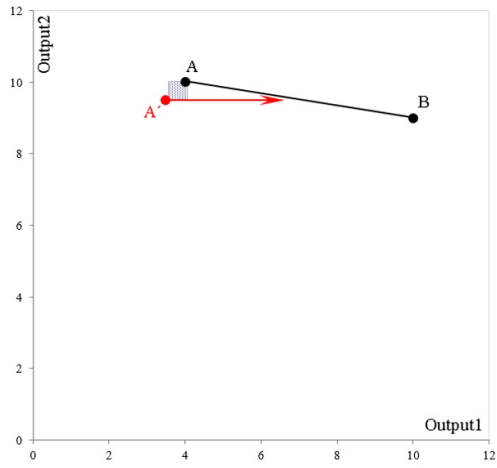
Figure 5: Benchmark $A'$ (a).

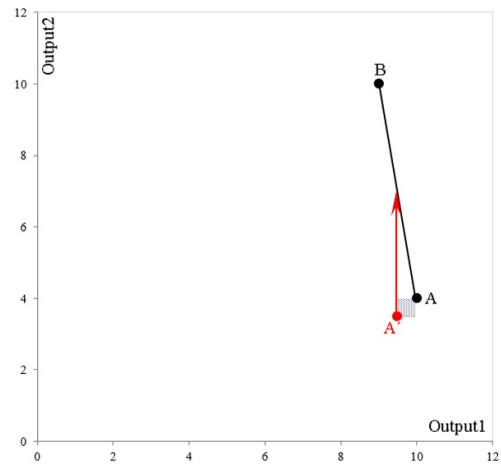
Figure 6: Benchmark $A'$ (b).

In other words, the selected neighbor of A is suggested to B, instead of suggesting to A. This phenomenon is due to the fact that the situation of B is better than the situation of A. Indeed, the outputs components of A are 4 and 10, while the corresponding outputs components of B are 10 and 9, respectively. The second outputs of A and B are almost the same, whereas the first output of A is 4 and the first output of B is 10. Therefore, when a neighbor of A such as $A'$ is selected, it is not in a corner of the frontier and it can be moved toward the best situation in the PPS. This is the Arash Method (AM) to look at the situation of an evaluated DMU and its neighbors at the same time, and measure the real relative efficiency of the DMU regardless whether it is in a corner of the frontier or not. In this example, $A'$ in Figures 5 and 6 is suggested to the point (1; 7, 9.5) and (1; 9.5, 7), respectively. If the introduced errors are supposed as $\varepsilon_1^+ = 0.1$ and $\varepsilon_2^+ = 0.1$ or smaller, the same illustration can also be demonstrated, and such introduced neighbors of A are still benchmarked toward B. Moreover, the length of the line segment $AA'$ is almost 0.71 where $\varepsilon_1^+ = \varepsilon_2^+ = 0.5$, and it is almost 0.14 where $\varepsilon_1^+ = \varepsilon_2^+ = 0.1$. As can be seen, when the value of epsilon is smaller, the diameter of depicted square in Figures 5 and 6 is smaller, therefore the sum of slacks is also smaller, but the benchmark is still toward B. For more illustration about KAM see Khezrimotlagh *et al.* (2012a-b) and Khezrimotlagh (2014a).

Now, if the place of B is changed to $B_1$ (Figures 7 and 8) and $B_2$ (Figures 9 and 10), that is, the situation of B becomes weaker and weaker in comparison with the situation of A, the optimum slacks for $A'$ become smaller and smaller, respectively, as the arrows indicate.

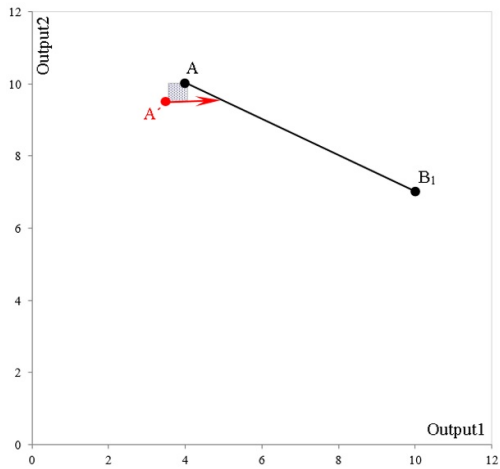
Figure 7: Benchmark $A'$ (c).

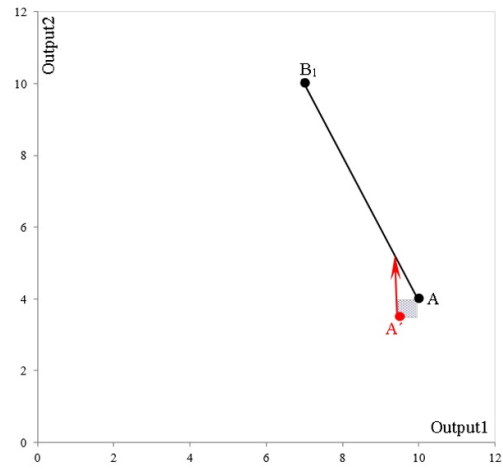
Figure 8: Benchmark $A'$ (d).

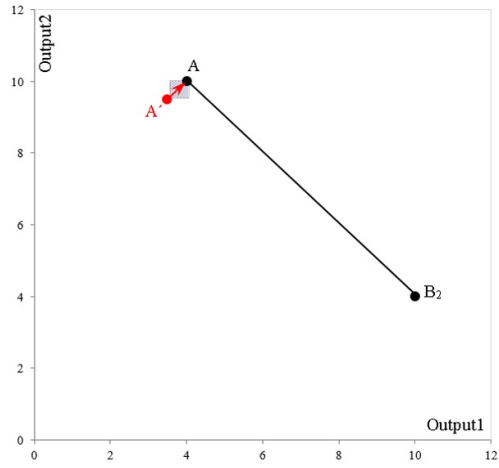
Figure 9: Benchmark $A'$ (e).

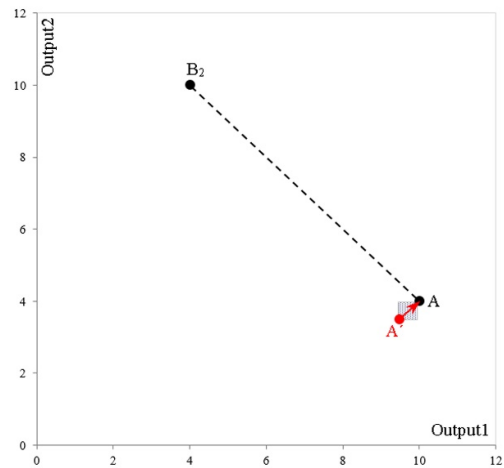
Figure 10: Benchmark $A'$ (f).

When the situation of B is changed to $B_2$, $A'$ is suggested to A, which shows that the situation of A is good, that is, the situation of A should be the same as the situation of $B_2$ or better. In this example, every point in the line segment $AB_2$ has the same efficiency score. So, there is no difference between efficiency score of A and $B_2$ while KAM is applied. When KAM efficiency scores of two technically efficient DMUs are the same with an introduced degree of freedom or error in inputs and outputs, it means that the situation of both DMUs are the same on the frontier while that error is introduced.

For every two sets of DMUs the above simple method can be exemplified. If in a practice, the frontier is the same as $AB_2$, it is obvious that there is neither difference between the efficiency scores of A and $B_2$ nor other points on the line segment $AB_2$, and as a result the neighbors of A or $B_2$ are suggested to A or $B_2$, respectively. However, if there is a very small difference between the efficiency score of two DMUs such as AB or $AB_1$, the neighbors of A are benchmarked toward B or $B_1$. From this real phenomenon, KAM can discriminate between DMUs appropriately, and rectify the problem of all current DEA methods to arrange DMUs. Figures 11 and 12 depict all the above illustrations simultaneously.

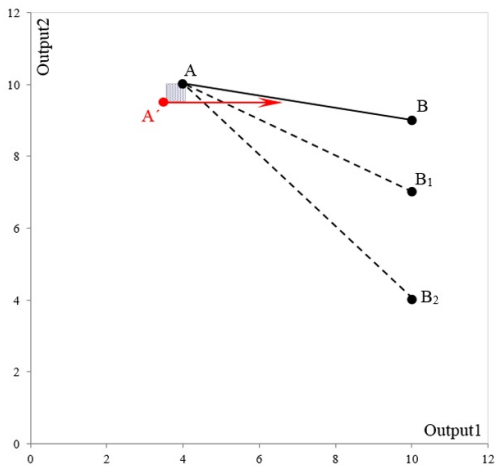
Figure 11: Different frontiers in output space (a).

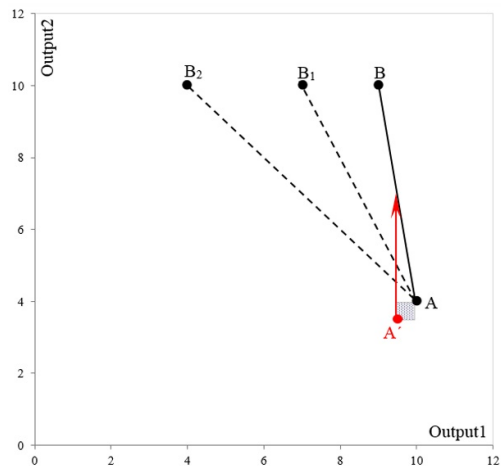
Figure 12: Different frontiers in output space (b).

Now, suppose that the number of DMUs is increased, as shown in Figure 13. The above illustrated method is easily able to distinguish between DMUs A-H. For instance, A has less efficiency score than B, B has less efficiency score than C, and C has less efficiency score than D and so on.

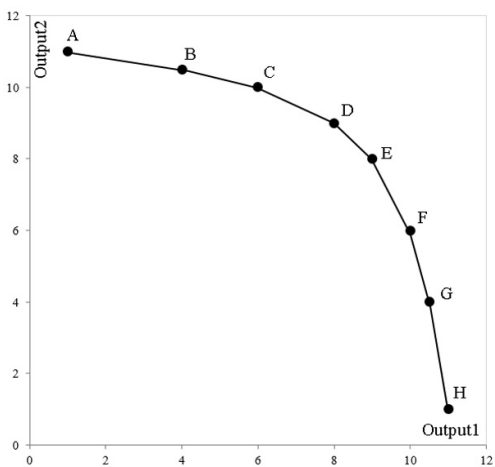
Figure 13: Eight technically efficient DMUs (a).

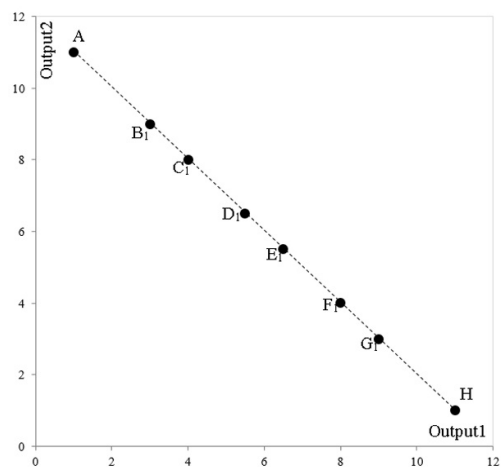
Figure 14: Eight technically efficient DMUs (b).

While a technically efficient DMU is under evaluation, KAM considers one of its neighbors according to the goal of practice, and checks where that neighbor (and of course that technically efficient DMU when the value of epsilon is negligible) should be benchmarked. KAM compares the efficiency score of the DMU and the efficiency score of the target of its close neighbor with a small negligible epsilon error, and defines a robust score for the DMU with epsilon degree of freedom.

If increasing the number of DMUs causes a situation such as DMUs in Figure 14, there are no differences between DMUs and the KAM efficiency score is the same for each DMU. However, if DMUs in Figure 15 are selected, the differences between A, B2-G2 and H are increased. All of these different situations can easily be measured by KAM without any computational complexities or additional methodologies. Figure 16 compares the different frontier while the number of technically efficient DMUs are increased. These different frontiers are depicted to understand how to describe the scores of KAM.

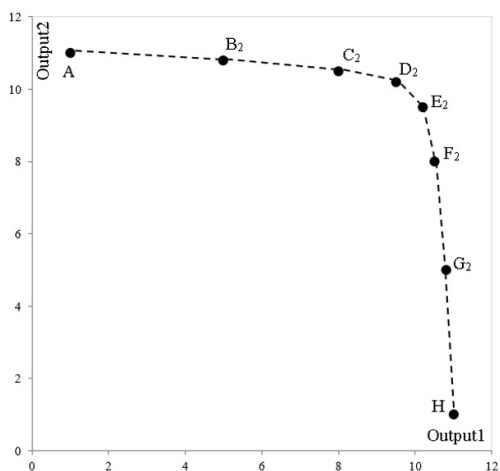
Figure 15: Eight technically efficient DMUs (c).

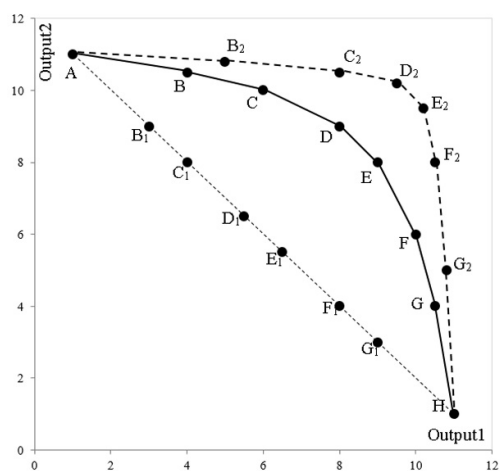
Figure 16: Eight technically efficient DMUs (d).

Note that, in the above illustration, the technically efficient DMUs are only considered, because when there are no sufficient numbers of DMUs in comparison with the numbers of factors, the numbers of technically efficient DMUs are increased. However, this is not a concern to discriminate between DMUs while KAM is applied.

Moreover, if the differences between the efficiency scores of a DMU and its neighbors are great, the DMU should not be known as an efficient DMU, even if it is a technically efficient DMU (See Definition in Section 2). Khezrimotlagh *et al.* (2013a) quoted that the technical efficiency is a necessary condition for being efficient and it is not enough to call a DMU "fully efficient". Unfortunately, the literature of DEA is full of incorrect interpretations considering technical efficiency as efficiency.

For instance, redial models are rarely suitable for measuring the relative efficiency scores, benchmarking and ranking DMUs. However, they have been used for almost all DEA theoretical and applicable studies. The redial approaches suppose that all points on the frontier have the same efficiency scores, which is absolutely not correct. This is a sad part of the history of DEA, if the concepts of super-efficiency models are remembered. Indeed, in one view, it is supposed that the points on the frontier have the same efficiency scores, so inefficient DMUs can be benchmarked to any of the points on the frontier. However, in another view, it is supposed that there are some differences between the points on the frontier, so they have to be ranked (the concept of super-efficiency models).

If the points on the frontier really have the same scores, why should one find different ranks for them? This transparent contradiction clearly shows the confusion in the literature of DEA, due to not considering the simple knowledge of logic. The ironic part of this argument is that DEA is judging whether factories do their job right, while DEA is not formulated utilizing simple logic. Moreover, as Khezrimotlagh *et al.* (2012a, 2013a, d) proved, even the super-efficiency DEA models are not valid to distinguish between technically efficient DMUs.

In short, KAM says that if a DMU has a good place on the frontier in comparison with other available DMUs, its neighbors should be benchmarked to that DMU and the differences of their efficiency scores should approximately be the same. Similar to the confident interval in Statistics, an introduced epsilon error allows users to identify efficient DMUs among the technically efficient ones. KAM fairly measures the efficiency scores of DMUs, and identifies the most efficient DMU with reasonable ranking and benchmarking for each, whether the number of DMUs are sufficient or not.

The next section depicts the advantages of the above technique while the number of DMUs is 32 and the number of variables is 45.

## 4. A real-life numerical example

Let's consider the example of Osman *et al.* (2011) with 32 nurses and 45 factors. In order to decrease the number of factors, they classified 20 factors into 6 input groups and the remaining factors into 9 output groups by calculating the average of factors and rounding data with two decimal digits. It is clear that this type of assessment misses some of the information. Indeed, the 45 factors are decreased to 15 factors by measuring the average and rounding the data. However, using the technique of KAM is easily able to discriminate between these nurses without missing any of information.

Applying 0-KAM in Variable Returns to Scale (VRS) technology shows that all nurses were technically efficient except the second nurse who was inefficient. Note that, VRS is considered due to increase the number of technically efficient DMUs and to illustrate the robustness of KAM.

In order to apply KAM, suppose that $\varepsilon = 0.001$, that is, only one thousandth errors in each factor, which is quite negligible according to data in the practice of Osman *et al.* (2011). Select the components of epsilon vector as $\varepsilon_j^- = \varepsilon \times x_{lj}$, and $\varepsilon_k^+ = \varepsilon \times x_{lk}$, for $j = 1,2,...,20$ and $k = 1,2,...,25$, which means KAM measures the efficiency score of an evaluated DMU's neighbor which has only one tenth percentage errors in its data in comparison with the evaluated DMU. For instance, the distance of the selected neighbor from the first nurse in the PPS is only 0.022, which is clearly negligible, and therefore, the relative efficiency score of the first DMU should be approximately the same as the relative efficiency score of its selected neighbor. A technically efficient DMU with lower differences between its relative efficiency score (which is 1) and the relative efficiency score of its close neighbors, will get a higher rank.

Assume that the weights are considered as $w_j^- = 1/x_{lj}$, and $w_k^+ = 1/y_{lk}$, for $j = 1,2,...,20$ and $k = 1,2,...,25$.

Now, let's apply 0.001-KAM in VRS for these 32 nurses with 20 input and 25 output variables. The results are depicted in Figure 17, which was sorted from most to least efficient nurses.

Nurses 32, 30, 25, 9, 8 and 22 are the first six efficient nurses and nurses 5, 11, 6, 12, 29 and 2 are the last six. The average of 32 efficiency scores is 0.996514 with the standard division of 0.016887, and the interquartile mean is 0.999657 with the standard division of 0.000105. According to Definition of Section 2, 0.001-KAM only knows the first four nurses as efficient with 0.001-DF where $\delta$ **= 0.0001**. It means, if only tenth percentage errors are introduced in data, only nurses 32, 30, 25 and 9 have good combinations of their factors and other nurses, although technically efficient, should improve their efficiency.

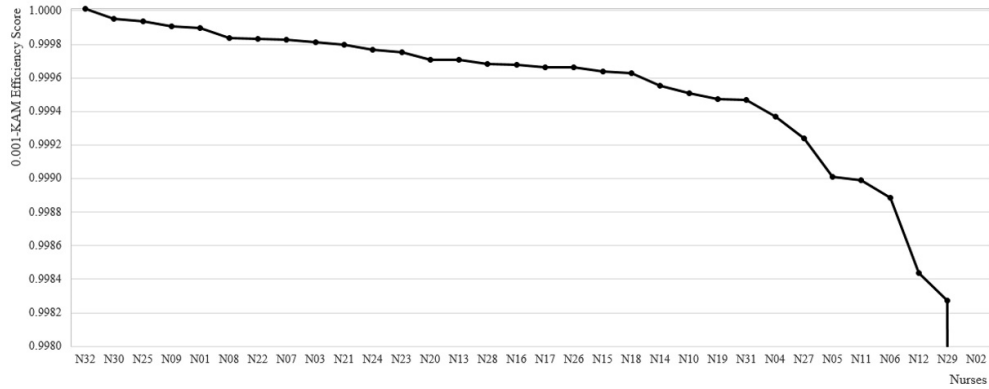
Figure 17: 0.001-KAM efficiency scores for 32 nurses with 45 variables.

For another experience, suppose that $\varepsilon$ **= 0.1** with the previous assumptions. Indeed, while a considered neighbor of an evaluated DMU is farther from the frontier, it can certainly be benchmarked to a better position on the frontier; however, the introduced degree of freedom shows the validity of considering the relative efficiency score of its neighbor as the relative efficiency score of that DMU.

Figure 18 depicts the results of 0.1-KAM which was sorted from the most to least efficient nurses. The average of efficiency scores is 0.956893 with standard division of 0.03565 and the interquartile mean is 0.967210 with the standard division of 0.009959. The first four nurses 32, 30, 25 and 9 still have the best rank, but nurse 1 is ranked higher than nurses 8 and 22, while 10 percentage errors are introduced in the data. Moreover, the technically efficient nurses 9 and 12 are ranked even lower than inefficient nurse 2 with 0.1-DF.

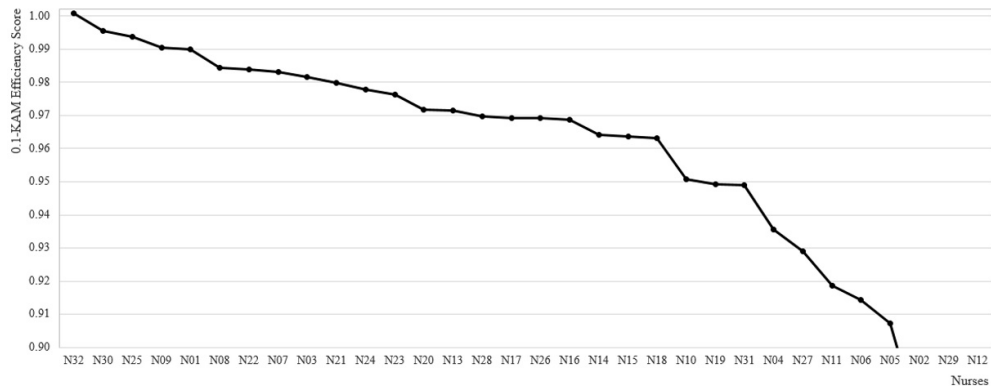
Figure 18: 0.1-KAM efficiency scores for 32 nurses with 45 variables.

These valuable outcomes are only one of the advantages of KAM to improve the discrimination power of DEA. KAM can consider a variety of weights for factors when the weights are unknown. If the weights of factors are available, the results of KAM is the same as the results of allocation models (Khezrimotlagh, 2014a). A discussion on how to select an epsilon for applying KAM can also be seen in (Khezrimotlagh 2014b).

## 5. Conclusion

This paper clearly illustrates how DEA is easily able to deal with a number of DMUs in comparison with a number of variables without any additional procedures, hybrid methodologies, computational complexities and extra conditions. The paper presents an obvious geometrical reason to prove that there are no concerns about the number of DMUs with an arbitrary number of variables in DEA while KAM is applied.